%% file: main.tex
\documentclass[12pt,a4paper,leqno]{article}

\usepackage{latexsym}
\usepackage{amssymb,exscale}
\usepackage[centertags]{amsmath}
\usepackage{amsthm}

\usepackage[dvips]{hyperref}

\input{preamble}

\begin{document}

\title{Moderate deviations for log-like functions\\
 of stationary Gaussian processes}

\author{Boris Tsirelson}

\date{}
\maketitle

\stepcounter{footnote}
\footnotetext{%
 This research was supported by \textsc{the israel science foundation}
 (grant No.~683/05).}

\begin{abstract}
A moderate deviation principle for nonlinear functions of Gaussian processes
is established. The nonlinear functions need not be locally
bounded. Especially, the logarithm is allowed. (Thus, small deviations of the
process are relevant.) Both discrete and continuous time is treated. An
integrable power-like decay of the correlation function is assumed.
\end{abstract}

\section*{Introduction}
\input{intro}

\section[Assumptions on the (nonlinear) function]
  {\raggedright Assumptions on the (nonlinear) function}
\label{sect1}
\input{sect1}

\section[Assumptions on the Gaussian process]
  {\raggedright Assumptions on the Gaussian process}
\label{sect2}
\input{sect2}

\section[The result]{\raggedright The result}
\label{sect3}
\input{sect3}

\section[Splitting the process]{\raggedright Splitting the process}
\label{sect4}
\input{sect4}

\section[A small deviation argument]{\raggedright A small deviation argument}
\label{sect5}
\input{sect5}

\section[Surgery]{\raggedright Surgery}
\label{sect6}
\input{sect6}

\section[Asymptotic variance]{\raggedright Asymptotic variance}
\label{sect7}
\input{sect7}

\section[Asymptotic exponential moments]{\raggedright Asymptotic exponential
 moments}
\label{sect8}
\input{sect8}

\bigskip
\filbreak
{
\small
\begin{sc}
\parindent=0pt\baselineskip=12pt
\parbox{4in}{
Boris Tsirelson\\
School of Mathematics\\
Tel Aviv University\\
Tel Aviv 69978, Israel
\smallskip
\par\quad\href{mailto:tsirel@post.tau.ac.il}{\tt
 mailto:tsirel@post.tau.ac.il}
\par\quad\href{http://www.tau.ac.il/~tsirel/}{\tt
 http://www.tau.ac.il/\textasciitilde tsirel/}
}

\end{sc}
}
\filbreak

\end{document}

%% file: preamble.tex
\numberwithin{equation}{section}

\swapnumbers
\theoremstyle{definition}
\newtheorem{theorem}[equation]{Theorem}
\newtheorem{lemma}[equation]{Lemma}
\newtheorem{example}[equation]{Example}

\renewcommand{\phi}{\varphi}

\newcommand{\I}{{\rm i}}
\newcommand{\D}{\mathrm{d}}
\newcommand{\E}{\mathrm{e}}

\renewcommand{\(}{\bigl(}
\renewcommand{\)}{\bigr)\vphantom{)}}
\renewcommand{\Pr}[1]{\mathbb{P}\mskip1.5mu\(\mskip1.5mu#1\mskip1.5mu\)}

\newcommand{\trace}{\operatorname{trace}}

\newcommand{\const}{\operatorname{const}}
\newcommand{\past}{\text{\textup{past}}}
\newcommand{\future}{\text{\textup{future}}}
\newcommand{\HS}{\text{\textup{HS}}}

\newcommand{\eps}{\varepsilon}
\newcommand{\si}{\sigma}
\newcommand{\ga}{\gamma}

\newcommand{\de}{\delta}
\newcommand{\al}{\alpha}
\newcommand{\be}{\beta}

\newcommand{\la}{\lambda}

\newcommand{\Ex}{\mathbb E\,}
\newcommand{\R}{\mathbb R}

\newcommand{\Z}{\mathbb Z}

\newcommand{\PR}[1]{\mathbb{P}\mskip1.5mu\bigg(\mskip1.5mu#1\mskip1.5mu\bigg)}

\newcommand{\neighborhood}[1]{$#1$\nobreakdash-\hspace{0pt}neighborhood}

\newcommand{\valued}[1]{$#1$\nobreakdash-\hspace{0pt}valued}
\newcommand{\periodic}[1]{$#1$\nobreakdash-\hspace{0pt}periodic}
\newcommand{\dimensional}[1]{$#1$\nobreakdash-\hspace{0pt}dimensional}

%% file: intro.tex
Questions on moderate deviations of random complex zeros \cite{So} lead
naturally to questions on moderate deviations of the logarithm of the absolute
value of a complex-valued Gaussian random field. Recently, Djellout, Guillin
and Wu established a moderate deviation principle for (nonlinear) functions of
dependent random variables (Gaussian, and more general)
\cite[Th.~2.7]{DGW}. However, their result does not answer the questions
mentioned above, for several reasons. The most important reason is that the
logarithm is not a differentiable function (nor even locally bounded).

The main result of the present work is another moderate deviation
principle. Unlike \cite{DGW}, I restrict myself to Gaussian processes, but
admit some non-locally-bounded functions (like the logarithm). My technique is
rather far from that of \cite{DGW}, and includes some arguments about small
deviations. Indeed, small values of the process lead to large negative values
of the logarithm.

%% file: sect1.tex
Let $ F : \R^d \to \R $ be a measurable function satisfying
\begin{gather}
\int F \, \D \ga^d = 0 \, , \label{1.1} \\
\int \E^F \, \D \ga^d < \infty \, , \quad \int \E^{-F} \, \D \ga^d < \infty \,
; \label{1.2}
\end{gather}
here and henceforth $ \ga^d $ is the standard Gaussian measure on $ \R^d $,
\[
\ga^d (\D x) = (2\pi)^{-d/2} \exp \( -\frac12 |x|^2 \) \, \D x \, .
\]
For any $ r \in (0,\infty) $ we define $ F_r : \R^d \to (-\infty,+\infty] $ by
\begin{equation}\label{1.25}
F_r (x) = \sup_{|y-x| \le r} F(y) \, .
\end{equation}
Here is the main assumption on $ F $: there exists $ C < \infty $ such that
\begin{equation}\label{1.3}
\int \exp \( F_r(x+0.5y) - F(x+0.5y) \) \, \ga^d (\D y) \le \E^{Cr}
\end{equation}
for all $ x \in \R^d $ and $ r \in (0,\infty) $.

Assumptions \eqref{1.1}, \eqref{1.2}, \eqref{1.3} will be referred to as
\emph{`the assumptions of Sect.~1'.}

\smallskip

Clearly, these assumptions are not invariant under replacement of $ F $ with $
aF $ for an arbitrary coefficient $ a \in (-\infty,0) \cup (0,\infty) $. (Even
$ a=-1 $ is not permitted by \eqref{1.3}.) However, our main results (Theorems
\ref{3.1}, \ref{3.2}) are evidently invariant under such replacement. Thus, we
could assume that $ aF $ satisfies \eqref{1.1}, \eqref{1.2}, \eqref{1.3}
for \emph{some} $ a \ne 0 $.

\begin{example}
Let $ F $ be a Lipschitz function, that is, $ |F(x)-F(y)| \le C|x-y| $ for all
$ x,y \in \R^d $. Then the function $ F(\cdot) - \int F \, \D \ga^d $
satisfies the conditions. Especially, \eqref{1.3} holds just because $ F_r(x)
\le F(x) + Cr $.
\end{example}

\begin{example}
Let $ d \ge 2 $. The function
\[
F(x) = \ln |x| - \int \ln |x| \, \ga^d(\D x)
\]
satisfies the conditions. Proof of \eqref{1.3}: if $ |y-x| \le r $ then $ F(y)
- F(x) = \ln \frac{|y|}{|x|} \le \ln \frac{|x|+|x-y|}{|x|} \le \ln \( 1 +
\frac{r}{|x|} \) $, thus, $ F_r(x) - F(x) \le \ln \( 1 + \frac{r}{|x|} \) $
and
\begin{multline*}
\int \exp \( F_r(0.5y) - F(0.5y) \) \, \ga^d (\D y) \le \int \Big( 1 +
\frac{2r}{|y|} \Big) \, \ga^d(\D y) = \\
= 1 + 2r \int \frac1{|y|} \, \ga^d (\D y) \le \exp \Big( 2r \int \frac1{|y|}
 \, \ga^d (\D y) \Big) \, ,
\end{multline*}
which is \eqref{1.3} for $ x=0 $. For arbitrary $ x $ we have
\[
\int \frac1{|x+y|} \, \ga^d (\D y) \le \int \frac1{|y|} \, \ga^d (\D y) \, ,
\]
which follows from the Anderson inequality, see for instance \cite[Th.~1.8.5
and Cor.~1.8.6]{Bo}.
\end{example}

Taking the limit $ r \to 0 $ in \eqref{1.3} we get $ \int | \nabla F (x+0.5y)
| \, \ga^d(\D y) \le C $ for a smooth $ F $; by approximation, \eqref{1.3}
implies that the first derivatives of $ F $ are locally finite measures. For
the one-dimensional case ($ d = 1 $) it means that \eqref{1.3} can hold only
for locally bounded $ F $. (However, $ F $ need not be continuous.)
Especially, the function $ x \mapsto \ln |x| $ on $ \R $ violates
\eqref{1.3}.

%% file: sect2.tex
\textsc{Discrete time}

\smallskip

Similarly to \cite{DGW} we consider a process $ (X_n)_{n\in\Z} $ that can be
written in the form
\begin{equation}\label{2.1}
X_n = \sum_{j\in\Z} a_{j-n} \xi_j = \sum_{j\in\Z} a_j \xi_{n+j}
\end{equation}
(`moving average process'), where $ \xi_n $ are independent. (About
convergence of the series, see below.) Unlike \cite{DGW} we assume that each $
\xi_n $ is an \valued{\R^d} random variable distributed $ \ga^d $, and each $
a_j $ is a matrix $ d \times d $. We assume that each $ X_n $ is also
distributed $ \ga^d $; it means that $ \sum a_j a^*_j $ is the unit
matrix. (Here $ a^* $ is the conjugate matrix.) The main assumption:
\begin{equation}\label{2.2}
a_j = O \Big( \frac1{|j|^{1.5+\eps}} \Big) \quad \text{for some } \eps>0 \, .
\end{equation}

The said above will be referred to as \emph{`the assumptions of Sect.~2 for
discrete time'.}

\smallskip

It follows from \eqref{2.2} that $ \sum \| a_j \| < \infty $, which is more
than enough for convergence of the series \eqref{2.1}. The Fourier transform
\[
g(\theta) = \sum_{n\in\Z} a_n \E^{\I n\theta}
\]
is a continuous \periodic{2\pi} matrix-valued function. The same holds for the
spectral density $ f $,
\[
f(\theta) = \frac1{2\pi} g(\theta) g^*(\theta) \, ; \quad \int_{-\pi}^\pi
\E^{\I n\theta} f(\theta) \D\theta = \Ex (X_n X_0^*) \, .
\]
It follows from \eqref{2.2} that $ \| g(\theta) - g(\eta) \| =
O(|\theta-\eta|^{0.5+\eps}) $ (see for instance \cite[Sect.~11.3]{Ba}). On the
other hand, every twice continuously differentiable (matrix-valued) function $
g $ is Fourier transform of a sequence $ (a_j)_j $ satisfying \eqref{2.2} for
$ \eps = 0.5 $.

If $ f(\cdot) $ is continuously differentiable twice and $ \det f(\cdot) $
does not vanish then $ f(\cdot) $ is of the form $ f(\theta) = (2\pi)^{-1}
g(\theta) g^*(\theta) $ with $ g(\cdot) $ satisfying \eqref{2.1}. (Just take
the positive square root of the positive matrix $ f(\theta) $.) Thus, a
process with such a spectral density belongs to our class, provided that $
\int_{-\pi}^\pi f(\theta) \, \D\theta $ is the unit matrix, and the process is
centered (zero-mean).

\medskip

\textsc{Continuous time}

\smallskip

Here we consider a process $ (X_t)_{t\in\R} $ that can be written in the form
\begin{equation}\label{2.3}
X_t = \int a_{s-t} \, \D w_s = \int a_s \, \D w_{t+s} \, ,
\end{equation}
where $ (w_s)_{s\in\R} $ is the standard \dimensional{d} Brownian motion
(two-sided; the past $ s \le 0 $ and the future $ s \ge 0 $ are independent,
and $ w_0 = 0 $), and $ s \mapsto a_s $ is a continuous matrix-valued function
on $ \R $. (The matrices are of size $ d \times d $.) We assume that each $
X_t $ is distributed $ \ga^d $; it means that $ \int a_s a^*_s \, \D s $ is
the unit matrix. The main assumption: there exists $ \eps > 0 $ such that the
function $ s \mapsto (|s|+1)^{1.5+\eps} a_s $ is bounded and H\"older
continuous; that is,
\begin{gather}
\sup_{s\in\R} (|s|+1)^{1.5+\eps} \| a_s \| < \infty \, , \label{2.4}
 \\
\sup_{s\in\R,\de\in(0,1)} \frac{ \| (|s+\de|+1)^{1.5+\eps} a_{s+\de} -
 (|s|+1)^{1.5+\eps} a_s \| }{ \de^\eps } < \infty \, . \label{2.5}
\end{gather}

The said above will be referred to as \emph{`the assumptions of Sect.~2 for
continuous time'.}

It follows from \eqref{2.4} that $ \int \| a_s \| \, \D s < \infty $, which is
more than enough for the linear stochastic integrals \eqref{2.3} to be
well-defined. The Fourier transform
\[
g(\la) = \int_{-\infty}^\infty a_s \E^{\I\la s} \, \D s
\]
is a continuous matrix-valued function. The same holds for the spectral
density $ f $,
\[
f(\la) = \frac1{2\pi} g(\la) g^*(\la) \, ; \quad \int_{-\infty}^\infty
\E^{\I\la t} f(\la) \D\la = \Ex (X_t X_0^*) \, .
\]
It follows from \eqref{2.4}, \eqref{2.5} that
\[
g(\la) = O(|\la|^{-\eps}) \, , \quad g(\la) - g(\mu) = O(|\la-\mu|^{0.5+\eps})
\, .
\]
On the other hand, every twice continuously differentiable (matrix-valued)
function $ g $ such that $ g(\la) $ and $ g''(\la) $ are $ O(|\la|^{-1-\eps})
$ (as $ |\la| \to \infty $) is Fourier transform of a function $ s \mapsto a_s
$ such that the functions $ s \mapsto a_s $ and $ s \mapsto s^2 a_s $ are
bounded and H\"older continuous, thus, \eqref{2.4}, \eqref{2.5} are satisfied.

%% file: sect3.tex
\begin{theorem}\label{3.1}
Let a function $ F : \R^d \to \R $ satisfy the assumptions of Sect.~1, and a
process $ (X_n)_{n\in\Z} $ satisfy the assumptions of Sect.~2 for discrete
time. Then

\noindent (a) the following limit exists:
\[
\si^2 = \lim_{n\to\infty} \frac1n \Ex \( ( F(X_1)+\dots+F(X_n) )^2 \) \, ;
\]

\noindent (b) if $ \si \ne 0 $ then for every $ \be \in (0,0.5) $ and $ \la \in \R $,
\[
\frac1{n^{1-2\be}} \ln \Ex \exp \bigg( \frac{\la}{\si n^\be} (
F(X_1)+\dots+F(X_n) ) \bigg) \to \frac{\la^2}2 \quad \text{as } n \to \infty
\, .
\]
\end{theorem}

\begin{theorem}\label{3.2}
Let a function $ F : \R^d \to \R $ satisfy the assumptions of Sect.~1, and a
process $ (X_t)_{t\in\R} $ satisfy the assumptions of Sect.~2 for continuous
time. Then

\noindent (a) the following limit exists:
\[
\si^2 = \lim_{t\to\infty} \frac1t \Ex \bigg( \int_0^t F(X_s) \, \D s \bigg)^2
\, ;
\]

\noindent (b) if $ \si \ne 0 $ then for every $ \be \in (0,0.5) $ and $ \la \in \R $,
\[
\frac1{t^{1-2\be}} \ln \Ex \exp \bigg( \frac{\la}{\si t^\be} \int_0^t F(X_s)
\, \D s \bigg) \to \frac{\la^2}2 \quad \text{as } t \to \infty \, .
\]
\end{theorem}

\medskip

It follows by the G\"artner-Ellis theorem (see \cite[Sect.~8]{Ell}) that for
every $ c \in [0,\infty) $,
\begin{gather*}
\frac1{n^{1-2\be}} \ln \Pr{ F(X_1)+\dots+F(X_n) < -c \si n^{1-\be} } \to
 -\frac{c^2}2 \, , \\
\frac1{n^{1-2\be}} \ln \Pr{ F(X_1)+\dots+F(X_n) > c \si n^{1-\be} } \to
 -\frac{c^2}2
\end{gather*}
as $ n \to \infty $ (discrete time), and
\begin{gather*}
\frac1{t^{1-2\be}} \ln \PR{ \int_0^t F(X_s) \, \D s < -c \si t^{1-\be} } \to
 -\frac{c^2}2 \, , \\
\frac1{t^{1-2\be}} \ln \PR{ \int_0^t F(X_s) \, \D s > c \si t^{1-\be} } \to
 -\frac{c^2}2
\end{gather*}
as $ t \to \infty $ (continuous time).

%% file: sect4.tex
\textsc{Discrete time}

\smallskip

Given a process $ (X_n)_{n\in\Z} $ satisfying the assumptions of Sect.~2 for
discrete time, and its representation \eqref{2.1}, we may split the process in
two independent processes,
\begin{gather*}
X_n = X_n^\past + X_n^\future \, , \\
X_n^\past = \sum_{j\le0} a_{j-n} \xi_j \, , \quad X_n^\future = \sum_{j>0}
a_{j-n} \xi_j \, .
\end{gather*}

\begin{lemma}\label{4.1}
There exists $ \eps > 0 $ such that
\[
\sup_{n\ge0} \( (n+1)^{1+\eps} | X_{-n}^\future | \) < \infty \quad \text{and}
\quad
\sup_{n>0} \( n^{1+\eps} | X_n^\past | \) < \infty \quad \text{a.s.}
\]
\end{lemma}

\begin{proof}
For $ n\ge0 $, using \eqref{2.2},
\[
\Ex | X_{-n}^\future |^2 = \sum_{j=1}^\infty \trace (a_{j+n} a^*_{j+n}) \le
\const \cdot \sum_{j=1}^\infty \frac1{ (j+n)^{3+2\eps} } = O(n^{-2-2\eps}) \,
.
\]
Thus, $ \sum_{n\ge0} \Pr{ | X_{-n}^\future | > n^{-1-0.5\eps} } < \infty $
(since $ X_{-n}^\future $ is Gaussian); the statement on $ X^\future $
follows. The statement on $ X^\past $ is similar.
\end{proof}

\medskip

\textsc{Continuous time}

\smallskip

Given a process $ (X_t)_{t\in\R} $ satisfying the assumptions of Sect.~2 for
continuous time, and its representation \eqref{2.3}, we may split the process
in two independent processes,
\begin{gather*}
X_t = X_t^\past + X_t^\future \, , \\
X_t^\past = \int_{-\infty}^0 a_{s-t} \, \D w_s \, , \quad X_t^\future =
\int_0^{\infty} a_{s-t} \, \D w_s \, .
\end{gather*}

\begin{lemma}\label{4.2}
There exists $ \eps > 0 $ such that
\[
\sup_{t\ge0} \( (t+1)^{1+\eps} | X_{-t}^\future | \) < \infty \quad \text{and}
\quad \sup_{t\ge0} \( (t+1)^{1+\eps} | X_t^\past | \) < \infty \quad
\text{a.s.}
\]
\end{lemma}

The proof, given afterwards, uses the following (quite general) lemma.

\begin{lemma}\label{4.3}
Let $ \phi : [0,\infty) \to \R $ and $ \eps>0 $ be such that the function $ x
\mapsto (x+1)^{0.5+\eps} \phi(x) $ is bounded and H\"older continuous, that
is,
\begin{gather*}
\sup_{x\ge0} \( (1+x)^{0.5+\eps} |\phi(x)| \) < \infty \, , \\
\sup_{x\ge0,0<\de<1} \frac{ | (1+x+\de)^{0.5+\eps} \phi(x+\de) -
(1+x)^{0.5+\eps} \phi(x) | }{ \de^\al } < \infty
\end{gather*}
for a given $ \al \in (0,1] $. Then there exists $ C < \infty $ such that for
all $ y \in [0,\infty) $ and $ \de \in [0,1] $
\[
\int_0^\infty | (1+y+\de)^\eps \phi(y+\de+x) - (1+y)^\eps \phi(y+x) |^2 \, \D
x \le C \de^{2\al} \, .
\]
\end{lemma}

\begin{proof}
Denote $ \psi(x) = (1+x)^{0.5+\eps} \phi(x) $. Choosing for every $ x $ some $
z_x \in [y,y+\de] $ we have
\begin{multline*}
(1+y+\de)^\eps \phi(y+\de+x) - (1+y)^\eps \phi(y+x) = \\
= \frac{ (1+y+\de)^\eps }{ (1+y+\de+x)^{0.5+\eps} } \psi(y+\de+x) - \frac{
(1+y)^\eps }{ (1+y+x)^{0.5+\eps} } \psi(y+x) = A_x + B_x + C_x \, ,
\end{multline*}
where
\begin{align*}
A_x &= \frac{ (1+y+\de)^\eps }{ (1+y+\de+x)^{0.5+\eps} } \( \psi(y+\de+x) -
 \psi(z_x+x) \) \, , \\
B_x &= \bigg( \frac{ (1+y+\de)^\eps }{ (1+y+\de+x)^{0.5+\eps} } - \frac{
 (1+y)^\eps }{ (1+y+x)^{0.5+\eps} } \bigg) \psi(z_x+x) \, , \\
C_x &= \frac{ (1+y)^\eps }{ (1+y+x)^{0.5+\eps} } \( \psi(z_x+x) - \psi(y+x) \)
 \, .
\end{align*}
Due to the triangle inequality in $ L_2(0,\infty) $ it is sufficient to choose
$ z_x $ (measurable in $ x $) such that $ \int_0^\infty A_x^2 \, \D x =
O(\de^{2\al}) $, $ \int_0^\infty B_x^2 \, \D x = O(\de^{2\al}) $ and $
\int_0^\infty C_x^2 \, \D x = O(\de^{2\al}) $. We have
\[
\int_0^\infty A_x^2 \, \D x = (1+y+\de)^{2\eps} \int_0^\infty O \(
(y+\de-z_x)^{2\al} \) \frac{ \D x }{ (1+y+\de+x)^{1+2\eps} } = O(\de^{2\al})
\, .
\]
Similarly, $ \int_0^\infty C_x^2 \, \D x = O(\de^{2\al}) $. These two
statements hold irrespective of the choice of $ z_x \in [y,y+\de] $. Now we
choose $ z_x $ such that
\begin{multline*}
\frac{ (1+y+\de)^\eps }{ (1+y+\de+x)^{0.5+\eps} } - \frac{ (1+y)^\eps }{
(1+y+x)^{0.5+\eps} } = \de \frac{\D}{\D z} \bigg|_{z=z_x} \frac{ (1+z)^\eps }{
 (1+z+x)^{0.5+\eps} } = \\
= \de \frac{ (1+z_x)^\eps }{ (1+z_x+x)^{0.5+\eps} } \bigg( \frac\eps{1+z_x} -
 \frac{ 0.5+\eps }{ 1+z_x+x } \bigg) \, .
\end{multline*}
The bracketed difference is evidently bounded; $ \psi(z_x+x) $ is also
bounded, thus,
\begin{gather*}
|B_x| \le \const \cdot \de \frac{ (1+z_x)^\eps }{ (1+z_x+x)^{0.5+\eps} } \le
 \const \cdot \de \frac{ (1+y+\de)^\eps }{ (1+y+x)^{0.5+\eps} } \, ; \\
\int_0^\infty B_x^2 \, \D x = O(\de^2) \cdot (1+y+\de)^{2\eps} \int_0^\infty
 \frac{ \D x }{ (1+y+x)^{1+2\eps} } = O(\de^2) = O(\de^{2\al}) \, .
\end{gather*}
\end{proof}

\begin{proof}[Proof of Lemma \textup{\ref{4.2}}]
First, for every $ t>0 $, using \eqref{2.4},
\begin{multline*}
\Ex | X_{-t}^\future |^2 = \int_0^\infty \trace (a_{s+t} a^*_{s+t}) \, \D s
 \le \\
\le \const \cdot \int_0^\infty \frac{\D s}{ (s+t+1)^{3+2\eps} } =
 O((t+1)^{-2-2\eps}) \, .
\end{multline*}
Second, we note that Lemma \ref{4.3} holds also for vector-valued (and
matrix-valued) functions, and apply it to the function $ \phi(t) = a_t $ ($ t
\ge 0 $), $ \eps $ in place of $ \al $ and $ 1+\eps $ in pace of $ \eps $
(recall \eqref{2.5}). We get
\[
\sup_{t\ge0} \int_0^\infty \| (t+1+\de)^{1+\eps} a_{t+\de+s} - (t+1)^{1+\eps}
a_{t+s} \|^2_{\HS} \, \D s = O(\de^{2\eps})
\]
for $ \de \le 1 $; here $ \| \cdot \|_{\HS} $ is the Hilbert-Schmidt norm, $
\| a \|^2_{\HS} = \trace(aa^*) $. Thus,
\begin{gather*}
\sup_{t\ge0} \Ex | (t+1+\de)^{1+\eps} X^\future_{-t-\de} - (t+1)^{1+\eps}
 X^\future_{-t} |^2 = O(\de^{2\eps}) \, , \\
\sup_{t\ge0} \Ex | (t+1)^{1+\eps} X^\future_{-t} |^2 < \infty \, .
\end{gather*}
It follows that the sample paths of the Gaussian process $ \( (t+1)^{1+\eps}
X^\future_{-t} \)_{t\ge0} $ are locally bounded (in fact, continuous). The
corresponding estimations are uniform, thus (see for instance
\cite[Th.~7.1.2]{Bo})
\[
\sup_{n\ge0} \Ex \sup_{t\in[n,n+1]} | (t+1)^{1+\eps} X^\future_{-t} | < \infty
\]
and moreover,
\[
\sup_{n\ge0} \PR{ \sup_{t\in[n,n+1]} | (t+1)^{1+\eps} X^\future_{-t} | > C }
\]
decays rapidly as $ C \to \infty $, namely, it is $ O(\E^{-\de C^2}) $ for
some $ \de > 0 $ (which can be obtained by Fernique's theorem, see for
instance \cite[Th.~2.8.5]{Bo}). By the Borel-Cantelli lemma,
\[
\sum_{n=0}^\infty \PR{ \max_{t\in[n,n+1]} \( (1+t)^{1+\eps} | X^\future_{-t} |
\) > n^{\eps/2} } < \infty \, .
\]
Therefore
\[
\sup_{t\ge0} \( (t+1)^{1+0.5 \eps} | X_{-t}^\future | \) < \infty \quad
\text{a.s.}
\]
The statement on $ X^\past $ is similar.
\end{proof}

%% file: sect5.tex
\textsc{Discrete time}

\smallskip

Let $ (X_n)_{n\in\Z} $ be a process satisfying the assumptions of Sect.~2 for
discrete time.

\begin{lemma}\label{5.1}
For every $ \si \in (0,1) $ there exist $ m \in \{1,2,\dots\} $ and an
\valued{\R^d} stationary Gaussian process $ (Y_k)_{k\in\Z} $ such that the two
processes
\[
(X_{mk})_{k\in\Z} \quad \text{and} \quad (Y_k+\si\xi_k)_{k\in\Z}
\]
are identically distributed; here $ \xi_k $ are independent \valued{\R^d}
random variables, each distributed $ \ga^d $, and the process $
(\xi_k)_{k\in\Z} $ is independent of the process $ (Y_k)_{k\in\Z}
$.\footnote{%
 $ mk $ in $ X_{mk} $ is just the product of $ m $ and $ k $.}
\end{lemma}

\begin{proof}
Here is a condition sufficient (and necessary, in fact) for existence of such
$ (Y_k)_k $ (for given $ m $ and $ \si $): the spectral density $ f_m $ of the
process $ (X_{mk})_k $ should exceed the spectral density of the process $
(\si \xi_k )_k $. That is, we need
\[
f_m(\theta) \ge \frac1{2\pi} \si^2 I \, ;
\]
here $ I $ is the unit matrix, and the inequality means that all the
eigenvalues of the Hermitian matrix $ f_m(\theta) $ lie on $ [\si^2,\infty)
$.

We have
\[
f_m(\theta) = \frac1m \bigg( f \Big( \frac{\theta}{m} \Big) + f \Big(
\frac{\theta+2\pi}{m} \Big) + \dots + f \Big( \frac{\theta+2\pi(m-1)}{m} \Big)
\bigg) \, ,
\]
the spectral density $ f $ of the given process $ (X_n)_n $ being a continuous
\periodic{2\pi} matrix-valued function on $ \R $ such that $ \int_{-\pi}^\pi
f(\theta) \, \D\theta = I $ (recall Sect.~\ref{sect2}). Therefore $
f_m(\theta) \to (2\pi)^{-1} I $ (as $ m \to \infty $) uniformly in $ \theta
$. It follows that $ f_m(\theta) \ge (2\pi)^{-1} \si^2 I $ for all $ \theta $,
if $ m $ is large enough.
\end{proof}

\medskip

\textsc{Continuous time}

\smallskip

Let $ (X_t)_{t\in\R} $ be a process satisfying the assumptions of Sect.~2 for
continuous time.

\begin{lemma}\label{5.2}
For every $ \si \in (0,1) $ there exist $ t \in (0,\infty) $ and a
discrete-time \valued{\R^d} stationary Gaussian process $ (Y_k)_{k\in\Z} $
such that the two discrete-time processes
\[
(X_{tk})_{k\in\Z} \quad \text{and} \quad (Y_k+\si\xi_k)_{k\in\Z}
\]
are identically distributed; here $ \xi_k $ are independent \valued{\R^d}
random variables, each distributed $ \ga^d $, and the process $
(\xi_k)_{k\in\Z} $ is independent of the process $ (Y_k)_{k\in\Z}
$.
\end{lemma}

\begin{proof}
Similarly to the proof of Lemma \ref{5.1} we consider the spectral density $
f_t(\cdot) $ of the process $ (X_{tk})_k $ and prove the inequality $
f_t(\theta) \ge (2\pi)^{-1} \si^2 I $ (which is sufficient).

We have
\[
f_t(\theta) = \frac1t \sum_{k\in\Z} f \Big( \frac{\theta+2\pi k}{t} \Big) \, ,
\]
the spectral density $ f $ of the given process $ (X_t)_t $ being an
integrable continuous matrix-valued function on $ \R $ such that $ f(\la) \ge
0 $ for all $ \la $, and $ \int_{-\infty}^\infty f(\la) \, \D\la = I $ (recall
Sect.~\ref{sect2}).

For $ M $ large enough,
\[
\int_{-M}^M f(\la) \, \D\la \ge (\si^2+\eps) I \quad \text{for some } \eps > 0
\, .
\]
For $ t $ large enough,
\[
f_t(\theta) \ge \frac1{2\pi} \int_{-M}^M f(\la) \, \D\la - \frac1{2\pi} \eps I
\ge \frac1{2\pi} \si^2 I \quad \text{for all } \theta \, .
\]
\end{proof}

%% file: sect6.tex
\textsc{Discrete time}

\smallskip

Let $ (X_n)_{n\in\Z} $ be a process satisfying the assumptions of Sect.~2 for
discrete time, and $ (Y_n)_{n\in\Z} $ its independent copy. We apply the split
of Sect.~\ref{sect4} to both:
\[
X_n = X_n^\past + X_n^\future \, , \quad Y_n = Y_n^\past + Y_n^\future \, .
\]
The four processes $ X^\past , X^\future, Y^\past, Y^\future $ are
independent. The processes $ X^\past $ and $ Y^\past $ are identically
distributed; symbolically, $ X^\past \sim Y^\past $. Also $ X^\future \sim
Y^\future $. Thus, we have four identically distributed stationary processes:
\[
X = X^\past + X^\future \sim X^\past + Y^\future \sim Y^\past + X^\future \sim
Y^\past + Y^\future = Y \, .
\]

Let $ F $ be a function satisfying the assumptions of Sect.~\ref{sect1}. We
introduce
\begin{gather*}
S_n^\future (X) = \sum_{k=1}^n F(X_k) \, , \\
S_n^\past (X) = \sum_{k=0}^{n-1} F(X_{-k}) = \sum_{k=1}^n F(X_{-n+k}) \, ,
\end{gather*}
denote by $ \mu_n $ the distribution of $ S_n^\future (X) $, symbolically $
S_n^\future (X) \sim \mu_n $, and observe that $ S_n^\past(X) \sim \mu_n $ and
$ S_m^\past(X) + S_n^\future(X) \sim \mu_{m+n} $. Further, we introduce
\begin{gather*}
D_n^\past = S_n^\past (X^\past+Y^\future) - S_n^\past (X^\past+X^\future) \, ,
 \\
D_n^\future = S_n^\future (Y^\past+X^\future) - S_n^\future
 (X^\past+X^\future)
\end{gather*}
and observe that
\begin{multline}\label{6.05}
S_m^\past(X) + S_n^\future(X) = \\
= S_m^\past(X^\past+Y^\future) + S_n^\future (Y^\past+X^\future) - D_m^\past -
 D_n^\future \, .
\end{multline}
This fact is instrumental to our purpose, since the two random variables $
S_m^\past(X^\past+Y^\future) $, $ S_n^\future (Y^\past+X^\future) $ are
independent, distributed $ \mu_m, \mu_n $ respectively, and the distribution
of their sum is close to $ \mu_{m+n} $ as far as $ D_m^\past + D_n^\future $
is relatively small.

\begin{lemma}\label{6.1}
There exists $ \eps > 0 $ such that
\[
\sup_{n>0} \Ex \exp \( \eps |D_n^\past| \) < \infty \, , \quad
\sup_{n>0} \Ex \exp \( \eps |D_n^\future| \) < \infty \, .
\]
\end{lemma}

\begin{proof}
It is sufficient to prove that
\[
\sup_{n>0} \Ex \exp \( \eps D_n^\future \) < \infty \, ,
\]
since the distribution of $ D_n^\future $ is symmetric (around $ 0 $), and the
assumptions of Sect.~\ref{sect2} are invariant under time
reversal. Equivalently, we may prove that
\begin{equation}\label{6.2}
\sup_{n>0} \Pr{ D_n^\future > C } = O(\E^{-\eps C})
\end{equation}
for some $ \eps > 0 $ and all $ C > 0 $.

By Lemma \ref{4.1}, $ \sup_{k>0} \( k^{1+\eps} |X_k^\past| \) < \infty $
a.s. The same holds for $ Y^\past $. We consider events
\[
A_u : \quad \sup_{k>0} \( k^{1+\eps} |X_k^\past-Y_k^\past| \) \le u \, .
\]
Fernique's theorem (mentioned in Sect.~\ref{sect4}) gives us $ \de > 0
$ such that
\begin{equation}\label{6.25}
\Pr{ A_u } \ge 1 - 2 \E^{-\de u^2} \quad \text{for } u \in (0,\infty) \, .
\end{equation}
We introduce $ Z_{u,k} = F_{uk^{-1-\eps}} (X_k) - F(X_k) \ge 0 $ (where $
F_{uk^{-1-\eps}} $ means $ F_r $ of \eqref{1.25} for $ r = uk^{-1-\eps} $) and
$ Z_u = \sum_{k>0} Z_{u,k} \in [0,\infty] $.
The intersection of $ A_u $ and the event $ D_n^\future > C $ is contained in
the event $ Z_u > C $, since $ |X_k^\past-Y_k^\past| \le uk^{-1-\eps} $
implies
\[
F ( Y_k^\past + X_k^\future ) - F ( X_k^\past + X_k^\future ) \le
F_{uk^{-1-\eps}} (X_k) - F (X_k) \, .
\]
Thus,
\begin{equation}\label{6.3}
\sup_{n>0} \Pr{ D_n^\future > C } \le \Pr{ Z_u > C } + 1 - \Pr{ A_u } \, .
\end{equation}
Lemma \ref{5.1} for $ \si = 0.5 $ gives us $ m \in \{1,2,\dots\} $ and $
(Y_k)_{k\in\Z} $ such that $ (X_{mk})_k \sim (Y_k+0.5\xi_k)_k $. We have
\begin{gather*}
Z_u = \sum_{k>0} Z_{u,k} = \sum_{j=1}^m \sum_{k\ge0} Z_{u,mk+j} \, ; \\
\exp ( m^{-1} Z_u ) \le \frac1m \sum_{j=1}^m \exp \bigg( \sum_{k\ge0}
 Z_{u,mk+j} \bigg) \, .
\end{gather*}
For each $ j \in \{1,\dots,m\} $,
$ \sum_{k\ge0} Z_{u,mk+j} = \sum_{k\ge0} \( F_{u(mk+j)^{-1-\eps}} (X_{mk+j}) -
F(X_{mk+j}) \) $ is distributed like $ \sum_{k\ge0} \( F_{u(mk+j)^{-1-\eps}} (
Y_k + 0.5 \xi_k ) - F ( Y_k + 0.5 \xi_k ) \) $.
For every (nonrandom) sequence $ (y_k)_k $, using \eqref{1.3} and denoting the
constant $ C $ of \eqref{1.3} by $ C_F $,
\begin{multline*}
\Ex \exp \Big( \sum_{k\ge0} \( F_{u(mk+j)^{-1-\eps}} ( y_k + 0.5 \xi_k ) - F (
 y_k + 0.5 \xi_k ) \) \Big) = \\
= \prod_{k\ge0} \Ex \exp \( F_{u(mk+j)^{-1-\eps}} ( y_k + 0.5 \xi_k ) - F (
 y_k + 0.5 \xi_k ) \) \le \\
\le \prod_{k\ge0} \exp \( C_F u (mk+j)^{-1-\eps} \) = \exp \Big( C_F u
 \sum_{k\ge0} (mk+j)^{-1-\eps} \Big) \le \exp (Bu) \, ,
\end{multline*}
where $ B = C_F \sum_{k\ge0} (mk+1)^{-1-\eps} < \infty $. Therefore
\begin{multline*}
\Ex \exp \Big( \sum_{k\ge0} Z_{u,mk+j} \Big) = \\
= \Ex \exp \Big( \sum_{k\ge0} \( F_{u(mk+j)^{-1-\eps}} ( Y_k + 0.5 \xi_k ) - F
 ( Y_k + 0.5 \xi_k ) \) \Big) \le \exp(Bu)
\end{multline*}
and $ \Ex \exp ( m^{-1} Z_u ) \le \E^{Bu} $, which implies
\[
\Pr{ Z_u > C } \le \exp ( Bu - m^{-1}C ) \, .
\]
We return to \eqref{6.3} and \eqref{6.25}:
\[
\sup_{n>0} \Pr{ D_n^\future > C } \le \exp ( Bu - m^{-1}C ) + 2 \E^{-\de u^2}
\]
for every $ u \in (0,\infty) $. Taking $ u = \sqrt C $ we get \eqref{6.2}.
\end{proof}

\medskip

\textsc{Continuous time}

\smallskip

Let $ (X_t)_{t\in\R} $ be a process satisfying the assumptions of Sect.~2 for
continuous time, and $ F $ a function satisfying the assumptions of
Sect.~\ref{sect1}. We proceed similarly to the discrete-time case: $ X =
X^\past + X^\future $, $ Y = Y^\past + Y^\future $ etc.; $ S_t^\future (X) =
\int_0^t F(X_s) \, \D s $, $ S_t^\past (X) = \int_{-t}^0 F(X_s) \, \D s $; $
D_t^\future = S_t^\future (Y^\past+X^\future) - S_t^\future
(X^\past+X^\future) $, and similarly $ D_t^\past $. We get
\begin{multline}
S_s^\past(X) + S_t^\future(X) = \\
= S_s^\past(X^\past+Y^\future) + S_t^\future (Y^\past+X^\future) - D_s^\past -
 D_t^\future \, .
\end {multline}

\begin{lemma}
There exists $ \eps > 0 $ such that
\[
\sup_{t>0} \Ex \exp \( \eps |D_t^\past| \) < \infty \, , \quad
\sup_{t>0} \Ex \exp \( \eps |D_t^\future| \) < \infty \, .
\]
\end{lemma}

The proof being quite similar to that of Lemma \ref{6.1}, I give a
sketch. Events
\[
A_u : \quad \sup_{t>0} \( (t+1)^{1+\eps} |X_t^\past-Y_t^\past| \) \le u
\]
satisfy $ \Pr{ A_u } \ge 1 - 2 \E^{-\de u^2} $ by Lemma \ref{4.2} and
Fernique's theorem. We introduce $ Z_{u,t} = F_{u(t+1)^{-1-\eps}} (X_t) -
F(X_t) \ge 0 $, $ Z_u = \int_0^\infty Z_{u,t} \, \D t \in [0,\infty] $ and get
\[
\sup_{t>0} \Pr{ D_t^\future > C } \le \Pr{ Z_u > C } + 2 \E^{-\de u^2} \, .
\]
Lemma \ref{5.2} for $ \si = 0.5 $ gives $ T $ and $ (Y_k)_k $ such that $
(X_{Tk})_k \sim (Y_k+0.5\xi_k)_k $. We have
\begin{gather*}
Z_u = \int_0^\infty Z_{u,t} \, \D t = \int_0^T \D s \sum_{k\ge0} Z_{u,Tk+s} \,
 ; \\
\exp ( T^{-1} Z_u ) \le \frac1T \int_0^T \D s \exp \bigg( \sum_{k\ge0}
Z_{u,Tk+s} \bigg) \, .
\end{gather*}
We use \eqref{1.3} as before:
\[
\Ex \exp \Big( \sum_{k\ge0} \( F_{u(Tk+s+1)^{-1-\eps}} ( y_k + 0.5 \xi_k ) - F
( y_k + 0.5 \xi_k ) \) \Big) \le \exp (Bu) \, ,
\]
where $ B = C_F \sum_{k\ge0} (Tk+1)^{-1-\eps} < \infty $. Thus, $ \Ex \exp
\Big( \sum_{k\ge0} Z_{u,Tk+s} \Big) \le \exp (Bu) $; $ \Ex \exp ( T^{-1} Z_u )
\le \E^{Bu} $; $ \Pr{ Z_u > C } \le \exp ( Bu - T^{-1}C ) $. Finally, $ u =
\sqrt C $ leads to $ \sup_{t>0} \Pr{ D_t^\future > C } = O(\E^{-\eps C}) $.

%% file: sect7.tex
\textsc{Discrete time}

\smallskip

Here we prove Item (a) of Theorem \ref{3.1}.

Let $ \mu_n $ be the distribution of $ F(X_1) + \dots + F(X_n) $. According to
Sect.~\ref{sect6}, $ \mu_{m+n} $ is close to the convolution $ \mu_m * \mu_n $
in the following sense. There exist random variables $ S_{m,n}, S'_m, S''_n $
such that
\begin{equation}\label{7.05}
\begin{gathered}
S_{m,n} \sim \mu_{m+n} \, , \quad S'_m \sim \mu_m \, , \quad S''_n \sim \mu_n
 \, , \\
S'_m \text{ and } S''_n \text{ are independent} \, , \\
\Ex \exp \( \eps | S_{m,n} - S'_m - S''_n | \) \le C
\end{gathered}
\end{equation}
for some $ \eps > 0 $, $ C < \infty $ not depending on $ m,n $. Namely, we may
take $ S_{m,n} = S_m^\past(X) + S_n^\future(X) $, $ S'_m =
S_m^\past(X^\past+Y^\future) $, $ S''_n = S_n^\future (Y^\past+X^\future) $
and note that
\begin{multline*}
\Ex \exp \( \eps | S_{m,n} - S'_m - S''_n | \) \le \Ex \exp \( \eps |
 D_m^\past | + \eps | D_n^\future | \) \le \\
\le \( \Ex \exp (2\eps |D_m^\past|) \)^{1/2} \( \Ex \exp (2\eps |D_n^\future|)
 \)^{1/2} \le C
\end{multline*}
by \eqref{6.05} and Lemma \ref{6.1}, if $ \eps $ is small enough and $ C $ is
large enough.

\begin{sloppypar}
Also,
\begin{equation}\label{7.1}
\int \E^{|x|} \, \mu_1(\D x) < \infty
\end{equation}
by \eqref{1.2}. In this section we need only second moments: $ \int x^2 \,
\mu_1(\D x) < \infty $ and
\[
\sup_{m,n} \Ex | S_{m,n} - S'_m - S''_n |^2 < \infty \, .
\]
Taking into account that the expectations vanish by \eqref{1.1}, we use
orthogonality and the triangle inequality in the space $ L_2 $ of random
variables:
\[
\sup_{m,n} \big| \| S_{m,n} \| - \sqrt{ \| S'_m \|^2 + \| S''_n \|^2 } \big| <
\infty \, .
\]
Thus, the numbers
\[
\si_n^2 = \int x^2 \, \mu_n(\D x) = \Ex \( F(X_1) + \dots + F(X_n) \)^2
\]
satisfy
\[
\sup_{m,n} \big| \si_{m+n} - \sqrt{ \si_m^2 + \si_n^2 } \big| < \infty \, .
\]
Existence of $ \lim_k ( 2^{-k/2} \si_{2^k} ) $ could be deduced readily, but
existence of $ \lim_n ( n^{-1/2} \si_n ) $ needs more effort. Here are two
quite general lemmas.
\end{sloppypar}

\begin{lemma}\label{7.2}
Let numbers $ a_1,a_2,\dots \in [0,\infty) $ and $ \eps > 0 $ satisfy
\[
a_{m+n} \le \sqrt{ a_m^2 + a_n^2 } + \eps
\]
for all $ m,n \in \{1,2,\dots\} $. Then
\[
a_n \le \Big( a_1 + \frac{\sqrt2}{\sqrt2-1} \eps \Big) \sqrt n
\]
for all $ n $.
\end{lemma}

\begin{proof}
For $ k = 0,1,2,\dots $ consider
\[
b_k = \max_{n\le2^k} \frac{ a_n }{ \sqrt n } \, .
\]
For each $ n \in \{ 1,2,3,\dots,2^k \} $
\[
\frac{ a_{2^k+n} }{ \sqrt{2^k+n} } \le \frac{ \sqrt{a_{2^k}^2+a_n^2}+\eps }{
  \sqrt{2^k+n} } \le \frac{ \sqrt{ 2^k b_k^2 + n b_k^2 } + \eps }{
  \sqrt{2^k+n} } \le b_k + \frac{\eps}{\sqrt{2^k}} \, ,
\]
therefore
\[
b_{k+1} \le b_k + 2^{-k/2} \eps \, ; \quad b_k \le b_0 +
\frac{\sqrt2}{\sqrt2-1} \eps \, .
\]
However, $ b_0 = a_1 $.
\end{proof}

Similarly,
\begin{equation}\label{7.3}
\text{if} \quad a_{m+n} \ge \sqrt{ a_m^2 + a_n^2 } - \eps \quad \text{then}
\quad a_n \ge \Big( a_1 - \frac{\sqrt2}{\sqrt2-1} \eps \Big) \sqrt n \, .
\end{equation}

\begin{lemma}\label{7.4}
Let numbers $ a_1,a_2,\dots \in [0,\infty) $ satisfy
\[
\sup_{m,n} \big| a_{m+n} - \sqrt{ a_m^2 + a_n^2 } \big| < \infty \, .
\]
Then there exists $ \lim_n (a_n/\sqrt n) \in [0,\infty) $.
\end{lemma}

\begin{proof}
Denote the given supremum by $ C $. For any $ k \in \{1,2,\dots\} $ we may
apply Lemma \ref{7.2} (together with \eqref{7.3}) to the sequence $
(a_k,a_{2k},\dots) $, obtaining
\[
| a_{nk} - a_k \sqrt n | \le \frac{\sqrt2}{\sqrt2-1} C \sqrt n \, ; \quad
\bigg| \frac{ a_{nk} }{ \sqrt{nk} } - \frac{ a_{k} }{ \sqrt{k} } \bigg| \le
\frac{\sqrt2}{\sqrt2-1} \frac{ C }{ \sqrt k } \, .
\]
All limiting points of the sequence $ (a_{nk}/\sqrt{nk})_n $ belong to the
\neighborhood{O(1/\sqrt k)} of the number $ a_k/\sqrt k $. The same holds for 
all limiting points of the sequence $ (a_{n}/\sqrt{n})_n $, since for $ \theta
\in \{ 0,1,\dots,k-1 \} $
\begin{gather*}
\left| \, a_{kn+\theta} - \sqrt{
 \smash{a_{kn}^2 + a_\theta^2} \mathstrut
 } \, \right| \le C \, ; \\
\left| \, \frac{ a_{kn+\theta} }{ \sqrt{kn+\theta} } - \sqrt{
  \smash{\frac{kn}{kn+\theta} \frac{a^2_{kn}}{kn} +
    \frac{a_\theta^2}{kn+\theta}} \vphantom{\frac A b}
 } \, \right| \le \frac{ C }{ \sqrt{kn+\theta}
 } \, ; \\
\left| \, \frac{ a_{kn+\theta} }{ \sqrt{kn+\theta} } - \sqrt{
  \smash{(1+o(1)) \frac{a^2_{kn}}{kn} + o(1)} \vphantom{\frac A b}
 } \, \right| \le o(1) \, .
\end{gather*}
Let $ x $ be a limiting point of $ (a_n/\sqrt n)_n $, then
\[
\Big| x - \frac{a_k}{\sqrt k} \Big| \le \frac{\sqrt2}{\sqrt2-1} \frac{ C }{
\sqrt k }
\]
for all $ k $. Thus, $ a_k/\sqrt k \to x $.
\end{proof}

It remains to apply Lemma \ref{7.4} to the sequence $ (\si_n)_n $.

\medskip

\textsc{Continuous time}

\smallskip

Item (a) of Theorem \ref{3.2} is verified similarly to that of Theorem
\ref{3.1}. We consider the distribution $ \mu_t $ of $ \int_0^t F(X_s) \, \D s
$ and note that $ \mu_{s+t} $ is close to $ \mu_s * \mu_t $ similarly to
\eqref{7.05}. Also,
\begin{equation}\label{7.*}
\Ex \exp \bigg( \frac1t \bigg| \int_0^t F(X_s) \, \D s \bigg| \bigg) \le \Ex
\exp |F(X_0)| < \infty \, .
\end{equation}
The numbers
\[
\si_t^2 = \int x^2 \, \mu_t(\D x) = \Ex \bigg( \int_0^t F(X_s) \, \D s
\bigg)^2
\]
defined for $ t \in [0,\infty) $ satisfy
\begin{gather*}
\sup_{s\in(0,t)} \si_s < \infty \quad \text{for } t < \infty \, , \\
\sup_{s,t} \left| \, \si_{s+t} - \sqrt{
 \smash{\si_s^2 + \si_t^2} \mathstrut
} \, \right| < \infty \, .
\end{gather*}

\begin{lemma}\label{7.9}
Let a function $ a : [0,\infty) \to [0,\infty) $ satisfy
\[
\sup_{s,t} \left| \, a(s+t) - \sqrt{
 \smash{a^2(s) + a^2(t)} \mathstrut
} \, \right| < \infty
\]
and be bounded on $ [0,t] $ for some (therefore, every) $ t > 0 $. Then there
exists $ \lim_{t\to\infty} (a(t)/\sqrt t) \in [0,\infty) $.
\end{lemma}

The proof is similar to that of Lemma \ref{7.4}. For any $ t \in (0,\infty) $
we apply Lemma \ref{7.2} (and \eqref{7.3}) to the sequence $ \(a(nt)\)_n $. A
limiting point $ x $ of the function $ s \mapsto a(s)/\sqrt s $ is also
a limiting point of the sequence $ \( a(nt) / \sqrt{nt} \)_n $ (boundedness of
$ a(\cdot) $ on $ [0,t] $ is used here), and we get $ \big| x -
\frac{a(t)}{\sqrt t} \big| \le \frac{\sqrt2}{\sqrt2-1} \frac{C}{\sqrt t} $.

It remains to apply Lemma \ref{7.9} to the function $ t \mapsto \si_t $.

%% file: sect8.tex
\textsc{Discrete time}

\smallskip

Here we prove Item (b) of Theorem \ref{3.1}.

Recall the numbers $ \si \in (0,\infty) $ from Item (a), $ \be \in (0,0.5) $
from Item (b), $ \eps > 0 $ from \eqref{7.05}. Denote $ S_{m,n} - S'_m - S''_n
$ from \eqref{7.05} by $ R_{m,n} $. Taking into account that $ \Ex R_{m,n} = 0
$ we get
\[
\Ex \exp (u R_{m,n}) \le \exp (C\si^2 u^2)
\]
for some $ C < \infty $ and all $ u \in [-\eps,\eps] $. We define functions $
\phi_n $ by
\[
\phi_n (\la) = \int \exp \Big( \frac{\la}{\si n^\be} x \Big) \, \mu_n(\D x) =
\Ex \exp \Big( \frac{\la}{\si n^\be} \( F(X_1) + \dots + F(X_n) \) \Big) \in
(0,\infty] .
\]

\begin{lemma}\label{8.1}
For all $ m,n $ and $ \la $ such that $ |\la| \le \eps \si (m+n)^{\be/2} $,
\begin{multline*}
\phi_m^p \bigg( \frac1p \Big( \frac{m}{m+n} \Big)^\be \la \bigg)
 \phi_n^p \bigg( \frac1p \Big( \frac{n}{m+n} \Big)^\be \la \bigg)
 \exp \bigg( \! -C \frac1p \frac{\la^2}{(m+n)^{3\be/2}} \bigg) \le \\[1mm]
\le \phi_{m+n} (\la) \le \\[1mm]
 \le \phi_m^{1/p} \bigg( p \Big( \frac{m}{m+n} \Big)^\be \la \bigg)
 \phi_n^{1/p} \bigg( p \Big( \frac{n}{m+n} \Big)^\be \la \bigg)
 \exp \bigg( C \frac{\la^2}{(m+n)^{3\be/2}} \bigg) \, ,
\end{multline*}
where
\begin{equation}\label{8.*}
p = \frac{ (m+n)^{\be/2} }{ (m+n)^{\be/2} - 1 } \, .
\end{equation}
\end{lemma}

\begin{proof}
The upper bound for $ \phi_{m+n} (\la) $: first,
\begin{multline*}
\phi_{m+n} (\la) = \Ex \exp \Big( \frac{\la}{\si (m+n)^\be} S_{m,n} \Big) =
 \Ex \exp \Big( \frac{\la}{\si (m+n)^\be} (S'_m+S''_n+R_{m,n}) \Big) \\
\le \Big\| \exp \Big( \frac{\la}{\si (m+n)^\be} (S'_m+S''_n) \Big)
 \Big\|_{L_p} \cdot \Big\| \exp \Big( \frac{\la}{\si (m+n)^\be} R_{m,n} \Big)
 \Big\|_{L_q} \, ,
\end{multline*}
where $ q = (m+n)^{\be/2} $. Second,
\begin{multline*}
\Big\| \exp \Big( \frac{\la}{\si (m+n)^\be} R_{m,n} \Big) \Big\|_{L_q} =
 \bigg( \Ex \exp \Big( q \frac{\la}{\si (m+n)^\be} R_{m,n} \Big) \bigg)^{1/q}
 \le \\
\le \exp \bigg( \frac1{(m+n)^{\be/2}} C\si^2 \Big( \frac{\la}{\si(m+n)^{\be/2}}
 \Big)^2 \bigg) = \exp \Big( C \frac{\la^2}{(m+n)^{3\be/2}} \Big) \, .
\end{multline*}
Third,
\begin{multline*}
\Big\| \exp \Big( \frac{\la}{\si (m+n)^\be} (S'_m+S''_n) \Big) \Big\|_{L_p}
 = \bigg( \Ex \exp \Big( p \frac{\la}{\si(m+n)^\be} (S'_m+S''_n) \Big)
 \bigg)^{1/p} = \\
= \bigg( \Ex \exp \Big( p \frac{\la}{\si(m+n)^\be} S'_m \Big) \bigg)^{1/p}
 \bigg( \Ex \exp \Big( p \frac{\la}{\si(m+n)^\be} S''_n \Big)
 \bigg)^{1/p} \, ,
\end{multline*}
and
\begin{multline*}
\Ex \exp \Big( p \frac{\la}{\si(m+n)^\be} S'_m \Big) = \Ex \exp \bigg(
 \frac1{\si m^\be} p \Big( \frac{m}{m+n} \Big)^\be \la S'_m \bigg) = \\
= \phi_m \bigg( p \Big( \frac{m}{m+n} \Big)^\be \la \bigg) \, ;
\end{multline*}
the same for $ S''_n $.

The lower bound for $ \phi_{m+n} (\la) $: first,
\begin{multline*}
\phi_m \bigg( \frac1p \Big( \frac{m}{m+n} \Big)^\be \la \bigg) \cdot
 \phi_n \bigg( \frac1p \Big( \frac{n}{m+n} \Big)^\be \la \bigg) = \\
= \Ex \exp \bigg( \frac1{\si m^\be} \frac1p \Big( \frac{m}{m+n} \Big)^\be \la
 S'_m \bigg) \cdot \Ex \exp \bigg( \frac1{\si n^\be} \frac1p \Big(
 \frac{n}{m+n} \Big)^\be \la S''_n \bigg) = \\
= \Ex \exp \bigg( \frac1p \frac{\la}{\si(m+n)^\be} (S'_m+S''_n) \bigg) = \Ex
 \exp \bigg( \frac1p \frac{\la}{\si(m+n)^\be} (S_{m,n} - R_{m,n}) \bigg) \le \\
\le \Big\| \exp \Big( \frac1p \frac{\la}{\si(m+n)^\be} S_{m,n} \Big)
 \Big\|_{L_p} \cdot \Big\| \exp \Big( -\frac1p \frac{\la}{\si(m+n)^\be}
 R_{m,n} \Big) \Big\|_{L_q} \, .
\end{multline*}
Second, the $ L_q $ norm is estimated by $ \exp \( C \frac1{p^2}
\frac{\la^2}{(m+n)^{3\be/2}} \) $ in the same way as before.
Third, the $ L_p $ norm is $ \( \Ex \exp \( \frac{\la}{\si(m+n)^\be} S_{m,n}
\) \)^{1/p} = \phi_{m+n}^{1/p} (\la) $. It remains to raise all that to the
power $ p $.
\end{proof}

Given a number $ \eps_1 \in (0,\eps] $, we consider (for every $ n $) the
smallest $ b_n $ and the largest $ a_n $ such that the inequalities
\begin{equation}\label{*}
\exp \( 0.5 n^{1-2\be} a_n \la^2 \) \le \phi_n(\la) \le \exp \( 0.5 n^{1-2\be}
b_n \la^2 \)
\end{equation}
hold for all $ \la $ satisfying $ |\la| \le \eps_1 \si n^{\be/2} $.

\begin{lemma}\label{8.2}
There exists $ N $ such that for all $ m,n $ satisfying $ N \le m \le n \le 2m
$,
\begin{align*}
a_{m+n} &\ge \frac{ (m+n)^{\be/2} - 1 }{ (m+n)^{\be/2} } \bigg( \frac{
 ma_m+na_n }{ m+n } - \frac{ 2C }{ (m+n)^{1-0.5\be} } \bigg) \, , \\
b_{m+n} &\le \frac{ (m+n)^{\be/2} }{ (m+n)^{\be/2} - 1 } \frac{
 mb_m+nb_n }{ m+n } + \frac{ 2C }{ (m+n)^{1-0.5\be} } \, .
\end{align*}
\end{lemma}

\begin{proof}
The bound for $ a_{m+n} $ will be verified for all $ m,n $. Let $ |\la| \le
\eps_1 \si (m+n)^{\be/2} $, then $ \big| \frac1p \(\frac{m}{m+n}\)^\be \la
\big| \le \eps_1 \si \frac{ m^\be }{ (m+n)^{\be/2} } \le \eps_1 \si m^{\be/2}
$ ($ p $ is still defined by \eqref{8.*}), thus
\begin{multline*}
\phi_m^p \bigg( \frac1p \Big( \frac{m}{m+n} \Big)^\be \la \bigg) \ge \exp
 \bigg( p \cdot \frac12 m^{1-2\be} a_m \cdot \frac1{p^2} \Big( \frac{m}{m+n}
 \Big)^{2\be} \la^2 \bigg) = \\
= \exp \bigg( \frac12 (m+n)^{1-2\be} \la^2 \cdot \frac1p \frac{ ma_m }{ m+n }
 \bigg) \, .
\end{multline*}
The same holds for $ \phi_n^p(\dots) $; Lemma \ref{8.1} gives
\[
\phi_{m+n} (\la) \ge \exp \bigg( \frac12 (m+n)^{1-2\be} \la^2 \cdot \frac1p
\Big( \frac{ ma_m+na_n }{ m+n } - \frac{ 2C }{ (m+n)^{1-0.5\be} } \Big) \bigg)
\, ,
\]
which verifies the bound for $ a_{m+n} $.

We take $ N $ such that $ 1.5^{-\be/2} + (2N)^{-\be/2} \le 1 $. Let $ N \le
m \le n \le 2m $, then
\[
\Big( \frac{n}{m+n} \Big)^{\be/2} \le \Big( \frac23 \Big)^{\be/2} \le 1 -
(2N)^{-\be/2} \le 1 - (m+n)^{-\be/2} = \frac1p \, ,
\]
therefore $ |\la| \le \eps_1 \si (m+n)^{\be/2} $ implies $ \big| p
\(\frac{n}{m+n}\)^\be \la \big| \le \(\frac{n}{m+n}\)^{\be/2} |\la| \le \eps_1
\si n^{\be/2} $ and
\begin{multline*}
\phi_n^{1/p} \bigg( p \Big( \frac{n}{m+n} \Big)^\be \la \bigg) \le \exp
 \bigg( \frac1p \cdot \frac12 n^{1-2\be} b_n \cdot p^2 \Big( \frac{n}{m+n}
 \Big)^{2\be} \la^2 \bigg) = \\
= \exp \bigg( \frac12 (m+n)^{1-2\be} \la^2 \cdot p \frac{ nb_n }{ m+n } \bigg)
 \, .
\end{multline*}
The same holds for $ \phi_m^{1/p}(\dots) $; Lemma \ref{8.1} gives
\[
\phi_{m+n} (\la) \le \exp \bigg( \frac12 (m+n)^{1-2\be} \la^2 \cdot \Big( p
\frac{ mb_m+nb_n }{ m+n } + \frac{ 2C }{ (m+n)^{1-0.5\be} } \Big) \bigg)
\, ,
\]
which verifies the bound for $ b_{m+n} $.
\end{proof}

Here are two quite general lemmas.

\begin{lemma}\label{8.3}
Let numbers $ B_1, B_2, \dots \in [0,\infty) $, $ \al \in (0,\infty) $ and $ r
\in [0,1] $ satisfy
\[
B_{m+n} \le \frac{ (m+n)^{\al} }{ (m+n)^{\al} - r } \frac{
 mB_m+nB_n }{ m+n } + \frac{ r }{ (m+n)^\al }
\]
for all $ m,n $ such that $ m \le n \le 2m $. Then
\[
\sup_n B_n \le (1+C_\al r) B_1 + C_\al r
\]
for some $ C_\al $ that depends on $ \al $ only.
\end{lemma}

\begin{proof}
We choose integers $ n_0 < n_1 < \dots $ such that $ \frac43 \le
\frac{n_{k+1}}{n_k} \le \frac32 $ for all $ k $, and $ n_0 = 2 $. We consider
\[
M_k = \max (B_1,B_2,\dots,B_{2n_k}) \, .
\]
For every integer $ n \in (2n_k,3n_k] $
\[
B_n = B_{n_k+(n-n_k)} \le \frac{n^\al}{n^\al-r} M_k + \frac{r}{n^\al} \, .
\]
Taking into account that $ 2n_{k+1} \le 3n_k $ we get
\[
M_{k+1} \le \frac{ (2n_k+1)^\al }{ (2n_k+1)^\al - r } M_k + \frac{r}{
(2n_k+1)^\al } \, .
\]
Introducing
\[
P_k = \prod_{j=0}^{k-1} \frac{ (2n_j+1)^\al }{ (2n_j+1)^\al - r } \, , \quad
P_0 = 1 \, ,
\]
we have $ P_k \ge 1 $, $ P_k \to P_\infty < \infty $, and $ M_{k+1} \le
\frac{P_{k+1}}{P_k} M_k + \frac{r}{ (2n_k+1)^\al } $, $
\frac{M_{k+1}}{P_{k+1}} \le \frac{M_k}{P_k} + \frac{r}{ (2n_k+1)^\al } $,
therefore
\begin{gather*}
\sup_k \frac{ M_k }{ P_k } \le \frac{ M_0 }{ P_0 } + r \sum_{k=0}^\infty
 (2n_k+1)^{-\al} \, ; \\
\sup_n B_n = \sup_k M_k \le P_\infty M_0 + r P_\infty \sum_k (2n_k+1)^{-\al}
\, .
\end{gather*}

We note that
\[
B_2 \le \frac{ 2^\al }{ 2^\al - r } B_1 + \frac{ r }{ 2^\al } \le (1+C_\al r)
B_1 + C_\al r \, ;
\]
here and henceforth $ C_\al $ is \emph{some} constant that depends on $ \al $
only, not necessarily the same in all occurrences. Similarly,
\[
B_3 \le (1+C_\al r) \max(B_1,B_2) + C_\al r \le (1+C_\al r) B_1 + C_\al r \, ,
\]
the same for $ B_4 $, and we get
\[
M_0 = \max(B_1,B_2,B_3,B_4) \le (1+C_\al r) B_1 + C_\al r \, .
\]
Finally, $ n_k \ge 2(4/3)^k $, thus
\begin{gather*}
\ln P_\infty = -\sum_{k=0}^\infty \ln \( 1 - r (2n_k+1)^{-\al} \) \le
 -\sum_{k=0}^\infty \ln \( 1 - r (4(4/3)^k+1)^{-\al} \) \le C_\al r ;
 \\
P_\infty \le \exp(C_\al r) \le 1 + (\E^{C_\al}-1) r \, .
\end{gather*}
\end{proof}

\begin{lemma}\label{8.4}
Let numbers $ A_1, A_2, \dots \in [0,\infty) $, $ \al \in (0,\infty) $ and $ r
\in [0,1] $ satisfy
\[
A_{m+n} \ge \frac{ (m+n)^{\al} - r }{ (m+n)^{\al} } \cdot \frac{
 mA_m+nA_n }{ m+n } - \frac{ r }{ (m+n)^\al }
\]
for all $ m,n $ such that $ m \le n \le 2m $. Then
\[
\inf_n A_n \ge (1-C_\al r) A_1 - C_\al r
\]
for some $ C_\al $ that depends on $ \al $ only.
\end{lemma}

\begin{proof}
Let $ n_k $ and $ P_k $ be as in the proof of Lemma \ref{8.3}. We introduce $
M_k = \min(A_1,\dots,A_{2n_k}) $ and get
\begin{gather*}
M_{k+1} \ge \underbrace{ \frac{ (2n_k+1)^\al - r }{ (2n_k+1)^\al } }_{
 P_k/P_{k+1} } M_k - \frac{r}{ (2n_k+1)^\al } \, ; \\
P_{k+1} M_{k+1} \ge P_k M_k - P_{k+1} r (2n_k+1)^{-\al} \, ; \\
\inf_k P_k M_k \ge P_0 M_0 - P_\infty r \sum_k (2n_k+1)^{-\al} \, ; \\
\inf_n A_n = \inf_k M_k \ge \frac{M_0}{P_\infty} - r \sum_k (2n_k+1)^{-\al}
\ge (1-C_\al r) A_1 - C_\al r
\end{gather*}
similarly to the proof of Lemma \ref{8.3}.
\end{proof}

Recall that $ a_n,b_n $ defined by \eqref{*} depend implicitly on $ \eps_1 $.

\begin{lemma}\label{8.6}
For every $ \de > 0 $ there exist $ \eps_1 \in (0,\eps] $ and $ k \in \{
1,2,\dots \} $ such that
\[
\inf_n a_{kn} \ge 1 - \de \quad \text{and} \quad \sup_n b_{kn} \le 1 + \de \,
.
\]
\end{lemma}

\begin{proof}
Lemma \ref{8.2} shows that Lemma \ref{8.3} may be applied to $ B_n = b_{kn} $,
$ \al = \be/2 $ and $ r = k^{-\be/2} \max(1,2C) $ provided that $ k $ exceeds
the number $ N $ of Lemma \ref{8.2}. Therefore $ \sup_n b_{kn} \le (1+\de) b_k
+ \de $ for all $ k $ large enough. Similarly (using Lemma \ref{8.4}), $
\inf_n a_{kn} \ge (1-\de) a_k - \de $ for all $ k $ large enough. Also,
\[
1-\de \le \frac{ \Ex ( F(X_1)+\dots+F(X_k) )^2 }{ k \si^2 } \le 1+\de
\]
for all $ k $ large enough. After choosing such $ k $ we choose $ \eps_1 $
such that
\[
\exp \( (1-\de) \cdot 0.5 \phi''_k(0) \la^2 \) \le \phi_k(\la) \le \exp \(
(1+\de) \cdot 0.5 \phi''_k(0) \la^2 \)
\]
for all $ \la $ satisfying $ |\la| \le \eps_1 \si k^{\be/2} $; this is
possible, since $ \phi_k(0) = 1 $ and
\[
\phi'_k(0) = \frac1{\si k^\be} \Ex ( F(X_1)+\dots+F(X_k) ) = 0 \, .
\]
Taking into account that
\[
\phi''_k(0) = \frac1{\si^2 k^{2\be}} \Ex ( F(X_1)+\dots+F(X_k) )^2 \in [
(1-\de) k^{1-2\be}, (1+\de) k^{1-2\be} ]
\]
we get
\[
\exp \( (1-\de)^2 \cdot 0.5 k^{1-2\be} \la^2 \) \le \phi_k(\la) \le \exp \(
(1+\de)^2 \cdot 0.5 k^{1-2\be} \la^2 \) \, ,
\]
which means that $ b_k \le (1+\de)^2 $ and $ a_k \ge (1-\de)^2 $. Finally, $
\inf_n a_{kn} \ge (1-\de)^3 - \de $ and $ \sup_n b_{kn} \ge (1+\de)^3 + \de
$.
\end{proof}

We see that
\[
(1-\de) \frac{\la^2}2 \le \frac1{(kn)^{1-2\be}} \ln \phi_{kn}(\la) \le (1+\de)
\frac{\la^2}2
\]
for $ |\la| \le \eps_1 \si (kn)^{\be/2} $. Now we consider $ \phi_{kn+\theta}(\la)
$ for $ \theta \in \{0,1,\dots,k-1\} $ and $ |\la| \le 0.5 \eps_1 \si
(kn+\theta)^{\be/2} $, assuming that $ kn+\theta $ is large enough (namely, exceeds $
2^{2/\be} $). Similarly to the proof of Lemma \ref{8.2} we use Lemma
\ref{8.1}. Taking into account that $ p = \frac{ (kn+\theta)^{\be/2} }{
(kn+\theta)^{\be/2} - 1 } \le 2 $ and
\[
\Big| p \Big( \frac{kn}{kn+\theta} \Big)^\be \la \Big| \le p \Big( \frac{kn}{kn+\theta}
\Big)^\be \cdot 0.5 \eps_1 \si (kn+\theta)^{\be/2} \le \eps_1 \si (kn)^{\be/2}
\]
we get
\begin{multline*}
\phi_{kn+\theta} (\la) \le \phi_\theta^{1/p} \bigg( p \Big(
 \frac{\theta}{kn+\theta} \Big)^\be \la \bigg) \phi_{kn}^{1/p} \bigg( p \Big(
 \frac{kn}{kn+\theta} \Big)^\be \la \bigg) \exp \bigg( C
 \frac{\la^2}{(kn+\theta)^{3\be/2}} \bigg) \\
= \phi_\theta^{1/p} \( o(1) \) \cdot \exp \Big( \frac1p (kn)^{1-2\be} (1+\de)
 \frac12 p^2 (1+o(1)) \la^2 \Big) \exp \( o(1) \) = \\
= \exp \Big( \frac12 (kn+\theta)^{1-2\be} \la^2 (1+\de+o(1)) + o(1) \Big)
\end{multline*}
for large $ kn+\theta $. Similarly,
\[
\phi_{kn+\theta} (\la) \ge \exp \Big( \frac12 (kn+\theta)^{1-2\be} \la^2 (1-\de+o(1)) +
o(1) \Big) \, ,
\]
which completes the proof of Theorem \ref{3.1}(b).

\medskip

\textsc{Continuous time}

\smallskip

We define
\[
\phi_t (\la) = \int \exp \Big( \frac{\la}{\si t^\be} x \Big) \, \mu_t(\D x) =
\Ex \exp \Big( \frac{\la}{\si t^\be} \int_0^t F(X_s) \, \D s \Big) \in
(0,\infty] \, ;
\]
similarly to Lemma \ref{8.1},
\begin{multline*}
\phi_s^p \bigg( \frac1p \Big( \frac{s}{s+t} \Big)^\be \la \bigg)
 \phi_t^p \bigg( \frac1p \Big( \frac{t}{s+t} \Big)^\be \la \bigg)
 \exp \bigg( \! -C \frac1p \frac{\la^2}{(s+t)^{3\be/2}} \bigg) \le \\[1mm]
\le \phi_{s+t} (\la) \le \\[1mm]
 \le \phi_s^{1/p} \bigg( p \Big( \frac{s}{s+t} \Big)^\be \la \bigg)
 \phi_t^{1/p} \bigg( p \Big( \frac{t}{s+t} \Big)^\be \la \bigg)
 \exp \bigg( C \frac{\la^2}{(s+t)^{3\be/2}} \bigg) \, ,
\end{multline*}
where $ p = \frac{ (s+t)^{\be/2} }{ (s+t)^{\be/2} - 1 } $.
Given $ \eps_1 \in (0,\eps] $, we consider (for every $ t $) the smallest $
b_t $ and the largest $ a_t $ such that the inequalities 
\[
\exp \( 0.5 t^{1-2\be} a_t \la^2 \) \le \phi_t(\la) \le \exp \( 0.5 t^{1-2\be}
b_t \la^2 \)
\]
hold for all $ \la $ satisfying $ |\la| \le \eps_1 \si t^{\be/2} $.
Similarly to Lemma \ref{8.2}, there exists $ T $ such that for all $ s,t $
satisfying $ T \le s \le t \le 2s $,
\begin{equation}\label{8.8}
\begin{aligned}
a_{s+t} &\ge \frac{ (s+t)^{\be/2} - 1 }{ (s+t)^{\be/2} } \bigg( \frac{
 sa_s+ta_t }{ s+t } - \frac{ 2C }{ (s+t)^{1-0.5\be} } \bigg) \, , \\
b_{s+t} &\le \frac{ (s+t)^{\be/2} }{ (s+t)^{\be/2} - 1 } \frac{
 sb_s+tb_t }{ s+t } + \frac{ 2C }{ (s+t)^{1-0.5\be} } \, .
\end{aligned}
\end{equation}

\begin{lemma}
For every $ \de > 0 $ there exist $ \eps_1 \in (0,\eps] $ and $ t \in
(0,\infty) $ such that
\[
\inf_n a_{tn} \ge 1 - \de \quad \text{and} \quad \sup_n b_{tn} \le 1 + \de \,
.
\]
\end{lemma}

\begin{proof}
By \eqref{8.8}, Lemma \ref{8.3} may be applied to $ B_n = b_{tn} $,
$ \al = \be/2 $ and $ r = t^{-\be/2} \max(1,2C) $ provided that $ t $ exceeds
$ T $ of \eqref{8.8}. The rest is completely similar to the proof of Lemma
\ref{8.6}.
\end{proof}

We see that
\[
(1-\de) \frac{\la^2}2 \le \frac1{(tn)^{1-2\be}} \ln \phi_{tn}(\la) \le (1+\de)
\frac{\la^2}2
\]
for $ |\la| \le \eps_1 \si (tn)^{\be/2} $. Now we consider $
\phi_{tn+\theta}(\la) $ for $ \theta \in [0,t] $ and $ |\la| \le 0.5 \eps_1
\si (tn+\theta)^{\be/2} $, assuming that $ tn+\theta $ is large enough
(namely, exceeds $ 2^{2/\be} $). We proceed similarly to the discrete case,
taking into account that the functions $ \phi_\theta(\cdot) $ are continuous
at $ 0 $ uniformly in $ \theta \in [0,t] $, which follows from \eqref{7.*} and
convexity of these functions.